\documentclass[12pt]{article}

\usepackage{amssymb}
\usepackage{color}
\usepackage{graphicx}
\usepackage{amsfonts}
\usepackage{amsmath,amsthm}

\def\R{{\mathbb R}}
\def\Z{{\mathbb Z}}

\def\La {{\Lambda}}
\def\si {{\sigma}}

\def\ga{{ \gamma}}
\def\Ga{{ \Gamma}}
\def\eps{{ \epsilon}}

\def\om{{ \omega}}

\def\C{{\cal C}}
\def\L{{\cal L}}

\def\clk{\C_k^{\ell_-}}

\def\clz{\C_0^{\ell_-}}

\def\tb{\tilde b}
\def\th{\tilde h}
\def\tl{\tilde\ell}

\def\a{\alpha}
\def\b{\beta}
\def\g{\gamma}
\def\k{\kappa}
\def\s{\sigma}
\def\t{\tau}

\def\={&=&}
\def\nn{\nonumber}

\newcommand{\dis}{\displaystyle}

\newtheorem{thm}{Theorem}
\newtheorem{prop}{Proposition}

\newtheorem{lem}{\indent Lemma}

\newtheorem{cor}{Corollary}

\begin{document}

\vskip.5cm
\title {Phase transitions in layered systems}

\author{Luiz Renato Fontes\footnote{Instituto de Matem\'atica e
Estat\'\i stica. Universidade de S\~ao Paulo, SP, Brazil. E-mail:
lrfontes@usp.br},
Domingos H. U. Marchetti\footnote{Instituto de F\'\i sica. Universidade de
S\~ao Paulo, SP, Brazil. Email: marchett@if.usp.br},
Immacolata Merola\footnote{DISIM, Universit\`a di L'Aquila, L'Aquila, Italy.
Email: immacolata.merola@univaq.it},\\ Errico Presutti\footnote{GSSI,
L'Aquila, Italy. Email: errico.presutti@gmail.com}, and
Maria Eulalia Vares\footnote{Instituto de Matem\'atica. Universidade
Federal do Rio de Janeiro, RJ, Brazil. Email: eulalia@im.ufrj.br}}
 \maketitle

\begin{abstract}
We consider the Ising model on  $\mathbb Z\times \mathbb Z$  where
on each horizontal line $\{(x,i), x\in \mathbb Z\}$, called  ``layer'',
the interaction
is given by  a ferromagnetic Kac potential with coupling strength $J_\ga(x,y)=\ga J(\ga(x-y))$,
where $J(\cdot)$ is smooth and has compact support; we
then add a nearest neighbor ferromagnetic
vertical
interaction of strength $\ga^{A}$,  where $A\ge 2$ is fixed,
and prove that for any $\beta$ larger than the mean field critical value there
is  a phase transition for all $\ga$ small enough.

\end{abstract}




\section{Introduction}
\label{sec;1}
We consider the Ising model on the lattice $\mathbb Z\times \mathbb Z$,
denoting by $(x,i)$ its points.
On each horizontal line $\{(x,i), x\in \mathbb Z\}$, called the $i$-th ``layer'',
the interaction
is given by  a ferromagnetic Kac potential
so that the  interaction between the spins at $(x,i)$ and $(y,i)$ is
  \begin{equation}
  \label{1.1}
- \frac 12 J_\ga(x,y) \si(x,i)\si(y,i),
   \end{equation}
where $J_\ga(x,y) = \ga J(\ga(x-y))$, and $J(\cdot)$ is a symmetric
smooth probability density on $\R$ with compact support. To fix the
notation we suppose $J(r)=0$ for $|r|\ge 1$.  We denote by $H_\ga^0$
the Hamiltonian with only the interactions \eqref{1.1} on each
layer, so that different layers do not interact with each other.

We fix the inverse temperature
$\beta>1$ (recalling that $\beta=1$ is the mean field critical value).
Since each layer
is independent of the others and one dimensional, the system
with Hamiltonian  $H_\ga^0$ does not have phase
transitions while its mean field
version (as derived by the Lebowitz-Penrose analysis by
taking first the thermodynamic limit and then letting
$\ga\to 0$) has a phase transition with infinitely many
extremal states, each one determined by
fixing on each layer a magnetization $\pm m_\beta$, $m_\beta>0$ the positive solution
of the mean field equation
  \begin{equation}
  \label{1.2}
 m_\beta = \tanh\{\beta m_\beta\}.
   \end{equation}
Purpose of this paper is to study what happens if we
put a ``very
small nearest neighbor vertical interaction''
  \begin{equation}
  \label{1.3}
- \eps \;\si(x,i)\si(x,i+1).
   \end{equation}
We take hereafter $\eps = \ga^{A}$, where $A\ge 2$ is fixed,
and call $H_\ga$ the Hamiltonian with both  interactions, i.e.\ the horizontal
one, \eqref{1.1},  and the vertical one,   \eqref{1.3} with $\eps = \ga^{A}$.
The Lebowitz-Penrose
limit is the same for $H^0_\ga$ and $H_\ga$, i.e.\ it is
not changed by the interaction \eqref{1.3}.  However the behavior of the system when $\ga>0$ is fixed (and suitably small) is
completely different.  Let $\La$ be a square in $\mathbb R^2$
and $\mu_{\ga,\La }^{\rm per}$ the  Gibbs measure with
Hamiltonian $H_\ga$ on $\La \cap (\mathbb Z\times \mathbb Z)$ with
periodic boundary conditions.  

\medskip
\begin{thm}
\label{thm1.1}
Fix $\beta>1$. There exists $\gamma_0>0$ (depending on $\beta$), so that for all $\ga\in (0,\gamma_0)$
  \begin{equation}
  \label{1.4}
{\rm weak} \lim_{\La\to \mathbb R^2}  \mu_{\ga,\La }^{\rm per} = \frac 12
\Big( \mu^+_\ga + \mu^-_\ga\Big)
   \end{equation}
where $\mu^{\pm}_\ga$ is the DLR measure obtained by taking the thermodynamic limit with plus, respectively minus,
boundary conditions.  Also $\mu^{+}_\ga \ne \mu^-_\ga$ for all such $\ga$. Furthermore, the expected values
of the spins converge as $\ga\to 0$ to their mean field values:
  \begin{equation}
  \label{1.5}
  \lim_{\ga\to 0}  \mu_{\ga}^{\pm} (\si(x,i))= \pm m_\beta.
   \end{equation}

\end{thm}

\medskip

The proof is given in the next sections
and it is obtained by establishing the validity of the Peierls bounds
for contours which are defined on each layer  following the coarse-grained
procedure in  \cite{presutti}.
The strategy for proving phase transitions in
$d\ge 2$ Ising systems with Kac potentials, as in \cite{CP,BZ,presutti},  is to
prove that for $\ga$ small enough the weight of a contour is well approximated by
the corresponding free energy excess of the associated Lebowitz-Penrose functional.
This does not work here because,  due to the smallness of the vertical
interaction \eqref{1.3}, the Lebowitz-Penrose functional does not penalize
phase changes between contiguous layers.  The analysis of the
interaction among layers is the main original part of the present paper and it is based
on the following idea.

The typical configurations for the Hamiltonian $H^0_\ga$ are made on each layer
by sequences of intervals where the empirical averages of the spins are alternatively close to $m_\beta$
and $-m_\beta$, the length of such intervals scales as $e^{c\ga^{-1}}$ ($c$ a positive constant),
as it was first observed in \cite{COP}. 
If this behavior were to persist after the vertical interaction \eqref{1.3} it would make the
interaction among intervals of different phase in contiguous layers
of the order $\eps e^{c\ga^{-1}}$; if $\eps$ is a power of $\ga$, as in Theorem \ref{thm1.1},
the Gibbs factor would depress such configurations and this is behind our proof of
the Peierls bounds for contours which describe a phase change between contiguous layers.

We hope our present results will help attacking the following problems which
arise naturally from the above considerations:

\begin{itemize}

\item  What happens in the thermodynamic limit to the Gibbs measure $\mu_{\ga,\La}^{+,-}$ defined
by putting plus boundary conditions on the layers $i\ge 0$ and minus boundary conditions
on the layers $i<0$ ?  Is the limit a Dobrushin state, maybe when the layers are $d>1$ dimensional ?

\item  Does the system still have a phase transition when $\beta=1$ (i.e.\ the mean field critical value)
and the vertical interaction \eqref{1.3} has strength $\eps>0$
independent of $\ga$ but arbitrarily small ?

\item Does Theorem \ref{thm1.1} extend to the case when on each layer line we have
a system of hard rods with attractive Kac pair potentials and a small
attractive  vertical interaction as in \eqref{1.3} ?  If the answer is positive
this would be an example where the original Kac proposal for the liquid-vapor phase transitions
can be carried through.

\end{itemize}

\noindent {\bf Comments.}
The idea of considering a Kac type interaction in each layer combined with a fixed  nearest
neighbor interaction in the vertical direction is by no means new. The reader is referred to a paper
by Kac and Helfand \cite{kac-helfand} in the early sixties. See also \cite {kac-thompson}. What seems new
to us is the consideration of the multiplicity of Gibbs measures for fixed (and very small) values of this
vertical interaction, beyond the Lebowitz-Penrose limit.


\vskip1cm

\setcounter{equation}{0}

\section{Contours}
\label{sec:2}

Following Chapter 9 in \cite{presutti} we
implement the program outlined in the Introduction
by a  coarse graining procedure.   For any $\ell \in \{ 2^n,\;n\in \mathbb Z\}$,  $i\in \mathbb Z$ and  $ k\in \mathbb Z$
we set:
     \begin{equation}
        \label{2.1}
C^{\ell,i}_{k \ell}= \Big\{ (x,i): k\ell \le x < (k+1)\ell\Big\},\;\;\; C^{\ell,i}_{x}=C^{\ell,i}_{k \ell}\; {\rm if}\;
(x,i) \in C^{\ell,i}_{k \ell}
    \end{equation}
and call $\mathcal D^{\ell,i}=\{C^{\ell,i}_{k \ell}, k\in \mathbb Z\}$,
$\mathcal D^{\ell}=\{\mathcal D^{\ell,i}, i\in \mathbb Z\}$.

We shall use  three basic parameters, two lengths $\ell_{\pm}$
and an accuracy $\zeta>0$ which all depend on $\ga$:
    \begin{equation}
    \label{2.2}
\ell_{\pm}= \ga^{-(1 \pm \alpha)},\quad \zeta=
\ga^a,\qquad 1\gg \alpha\gg a>0 \,
     \end{equation}
supposing for notational simplicity that $\ell_{\pm}
\in \{2^{n}, n\in \mathbb N_+\}$: this is a restriction on
$\ga$  and $\alpha$ which could be removed by taking
integer parts in \eqref{2.2}.  We  shortly call  $\ell_{\pm}$ intervals
the intervals which belongs to $\mathcal D^{\ell_{\pm}}$.

Define the empirical magnetization on the scale $\ell_-$ as
    \begin{equation}
    \label{2.3}
\si^{(\ell_-)}(x,i) := \frac{1}{\ell_{-}} \sum_{y: (y,i)\in
C^{\ell_-,i}_{x}} \si(y,i).
     \end{equation}
The random variables $\eta(x,i)$, $\theta(x,i)$ and $\Theta(x,i)$ are then defined as follows:

\medskip

\begin{itemize}

\item\;   $\eta(x,i)=\pm 1$ if  $\dis{ \big|  \si^{(\ell_-)}(x,i) \mp m_\beta
 \big|\le \zeta}$
 and $=0$ otherwise.

\item\;   $\theta(x,i)= 1$, [$=-1$], if
$\eta(y,i)= 1$, [$=-1]$, for all $(y,i)\in
C_x^{\ell_{+},i}$ and $=0$ otherwise.

\item\;    $\Theta(x,i)= 1$, [$=-1$], if
$\theta(x,i)= 1$, [$=-1$], on $ C_x^{\ell_{+},i}\cup
 C_{x'}^{\ell_{+},i} \cup  C_{x''}^{\ell_{+},i}$, where the latter
are the $\ell_+$ intervals immediately to the right and to the left
of $C_{x}^{\ell_{+},i}$ and $=0$ otherwise.

\end{itemize}

\medskip

\noindent
The phase of a site  $(x,i)$ is   ``plus'' if $\Theta(x,i)= \Theta(x,i\pm 1) =1$,
it is   ``minus''   if  $\Theta(x,i)= \Theta(x,i\pm 1) =-1$ and it
is ``undetermined'' otherwise.
Thus given a spin configuration $\si$ we have
a plus, a minus and an
undetermined region.
Calling ``connected'' $(x,i)$ and $(y,j)$
iff $|x-y|\le 1$, $|i-j|\le 1$
it then follows (recalling the definition of $\Theta$) that the plus and minus
regions are disconnected from each other by
the undetermined region.

We shall restrict in the sequel
to spin configurations
such that $\Theta=1$ outside of a compact
(the case
when $\Theta=-1$ 
can be recovered
via spin flip). Given such a  $\si$,
we call ``contours'' the pairs
$\Ga=({\rm sp}(\Ga), \eta_\Ga)$, where ${\rm sp}(\Ga)$ is a maximal connected component
of the undetermined region, called ``the spatial support of $\Ga$'',
and $\eta_\Ga$ is the restriction of $\eta$ to ${\rm sp}(\Ga)$,  called ``the specification of $\Ga$''.

\begin{figure}
  \centering
  \includegraphics[width=16cm]{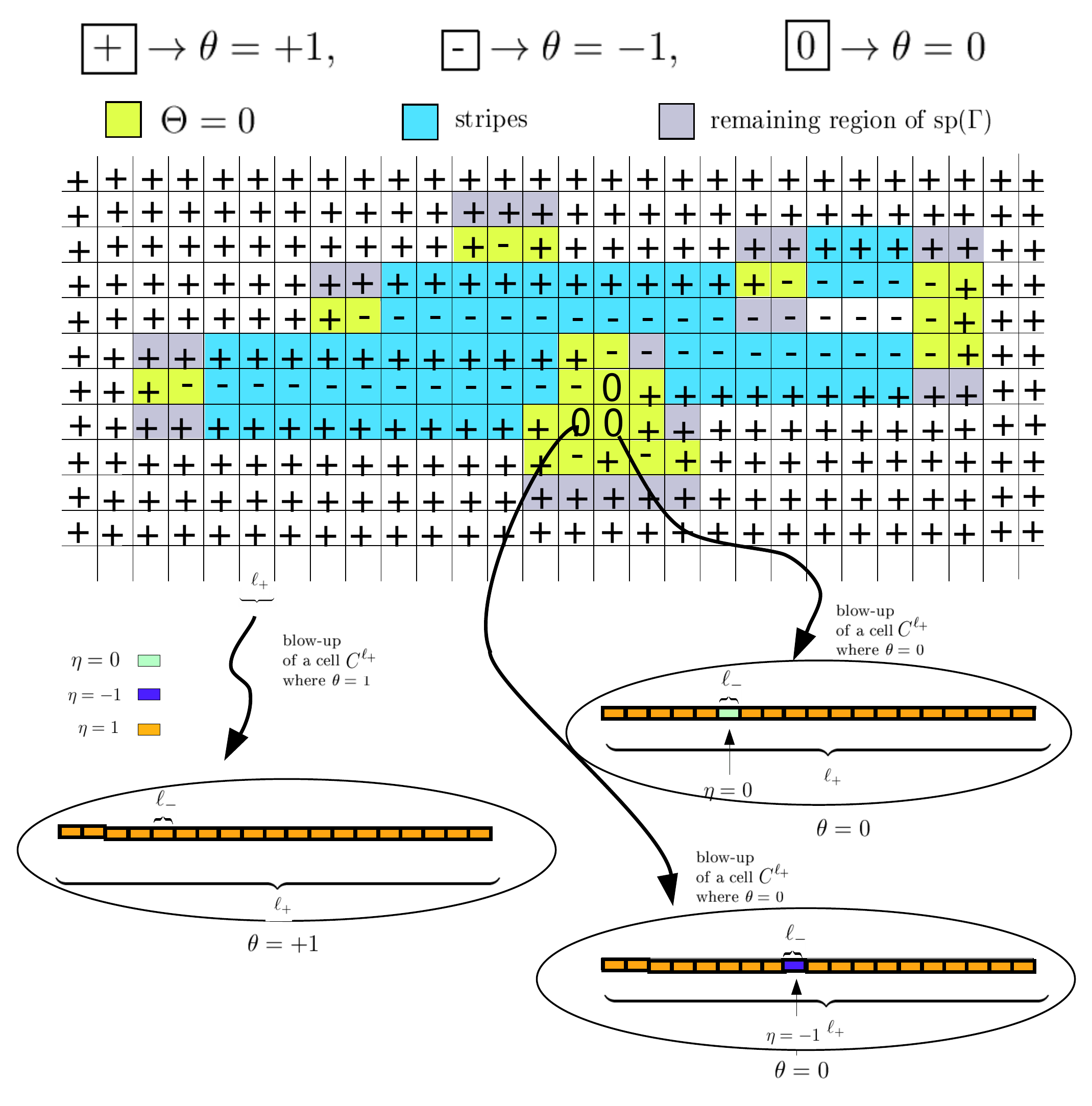}\\
  \caption{The picture illustrates the support of a contour, marked in the picture
by different colors. Each square represents an interval of length
$\ell_+$ on a line, and the values $+1,0,-1$ correspond to the
$\theta$ variables: the color yellow marks the intervals where
$\Theta=0$, the light-blue marks intervals where $\Theta$ is $+1$ or
$-1$ but in the line below or in the line above the sign is
opposite, i.e. the regions denoted as ``stripes" . The color grey
marks the intervals where $\Theta$ is $+1$ or $-1$ but in the line
below or in the line above $\Theta=0$. In the bottom we show some
blow-ups of these cells, where $\theta=1$ (the case $\theta=-1$ is
similar) and two examples where $\theta=0$.}\label{contour}
\end{figure}

Denote by ${\rm ext} (\Ga)$ the unbounded maximal connected
component of the complement of ${\rm sp}(\Ga)$ and $\partial_{\rm
ext}(\Ga)$ the union of all $\ell_+$ intervals in ${\rm ext} (\Ga)$
which are connected to ${\rm sp}(\Ga)$. Then (since  ${\rm sp}(\Ga)$
is bounded and connected) $\partial_{\rm ext}(\Ga)$ is connected;
moreover $\Theta \ne 0$ on $\partial_{\rm ext}(\Ga)$ (because ${\rm
sp}(\Ga)$ is a maximal connected component of the undetermined
region) and hence $\Theta$ is constant and different from 0 on
$\partial_{\rm ext}(\Ga)$ (because the plus and minus regions are
disconnected).  We shall call ``plus'' a contour $\Ga$ when
$\Theta=1$   on $\partial_{\rm ext}(\Ga)$ and ``minus'' otherwise.
(See Figure \ref{contour} for an illustration of a contour.)

Analogously we call ${\rm int}_k(\Ga)$ the bounded maximal connected components (if any)
of the complement of ${\rm sp}(\Ga)$, $\partial_{k}(\Ga)$ the $\ell_+$ intervals in
${\rm int}_k(\Ga)$ which are connected to ${\rm sp}(\Ga)$; then
$\Theta$ is constant and different from 0 on each $\partial_{k}(\Ga)$ and
we write $\partial^{\pm}_{k}(\Ga)$ if $\Theta=\pm 1$.  We also call
    \begin{equation}
    \label{2.4}
c(\Ga) = {\rm sp}(\Ga)\cup\bigcup_k {\rm int}_k(\Ga).
     \end{equation}
We are now ready to define the fundamental notion of ``weight of a contour''.  Let
$\Ga$ be a plus contour (the definition for minus contours is obtained by spin flip);
    \begin{equation}
    \label{2.5}
\{\si_{c(\Ga)} \Rightarrow \Ga \}= \Big\{ \si_{c(\Ga)}: \eta  = \eta_\Ga \;{\rm on}\;
 \, {\rm sp}(\Ga),\,\eta_\Ga = \pm 1 \;{\rm on \;all} \,\; \partial_{k}^{\pm}(\Ga)
 \Big\}
     \end{equation}
and $\bar \si_{\partial_{\rm ext}(\Ga)}$ a configuration such that $\Theta=1$ on the whole
$\partial_{\rm ext}(\Ga)$.  We then define  the weight of $\Ga$ with boundary conditions
$\bar \si_{\partial_{\rm ext}(\Ga)}$ as
    \begin{equation}
    \label{2.6}
W_\Ga(\bar\si_{\partial_{\rm ext}(\Ga)}):= \frac
{Z_{c(\Ga);\bar\si_{\partial_{\rm ext}(\Ga)}}(\{\si_{c(\Ga)} \Rightarrow \Ga \})}
{Z_{c(\Ga);\bar\si_{\partial_{\rm ext}(\Ga)}}(\{\si_{c(\Ga)}: \Theta = 1\; \rm{on}\;
{\rm sp}(\Ga)\; {and \; all}\, \partial_{k}^{\pm}(\Ga)\})},
     \end{equation}
where $Z_{\La,\bar \si_{\partial_{\rm ext}(\La)}}(C)$ is the partition function
in $\La$ with boundary conditions $\bar\si_{\partial_{\rm ext}(\La)}$ and constraint $C$.

The Peierls argument is based on (a) a bound on the weight of contours and (b) a counting argument
for the number of contours which contain a given site. These are  established in the next two sections.

\vskip1cm

\setcounter{equation}{0}

\section{Energy bounds}
\label{sec:3}

We shall prove here bounds on the weight of the contours which are exponentially small
with the exponent proportional to the spatial support $|{\rm sp}(\Ga)|$ of the contour.
To this end we introduce the notion of ``stripes'' in a contour $\Ga$.  The spatial support ${\rm sp}(S)$ of a stripe  $S$
is a set $\{(x,i): x \in I\}\cup \{(x,i+1): x \in I\}$ where $I= [k\ell_+, h\ell_+-1]$, $k,
h \in \mathbb Z$, $k<h$. $S$ is a $+-$ stripe in $\Ga$
if  ${\rm sp}(S) \subset {\rm sp}(\Ga) $ and:

\begin{itemize}

\item $\Theta = 1$ on the upper part of ${\rm sp}(S)$ and $=-1$ on the lower part ($\Theta$ is determined
on ${\rm sp}(\Ga)$ by the specification $\eta_\Ga$ of $\Ga$).

\item ${\rm sp}(S)$ is maximal with the above property, namely if $x \in [(k-1)\ell_+,k\ell_+)$
then at least one between $\Theta(x,i)$ and  $\Theta(x,i+1)$ is equal to 0 and the same holds for
$x \in [h\ell_+,(h+1)\ell_+)$.

\end{itemize}

$-+$ stripes are defined analogously (with $-$ on the top).  We call $|S|$ the number of sites
in the interval $I$ associated to the stripe $S$.  We have:

\medskip

\begin{thm}
\label{thm3.1} There is a positive constant $c$ so that for all
$\ga$ small enough the following holds. Let $\Ga$ be any plus
contour, $\mathcal S$ the set of all stripes in $\Ga$, $|\mathcal
S|$ the sum of $|S|$ over $S\in \mathcal S$ and $N_0$ the number of
intervals of $\mathcal D^{\ell_+}$ contained in ${\rm sp}(\Ga)$
where $\Theta=0$. Then for any $\bar\si_{\partial_{\rm ext}(\Ga)}$
such that $\Theta=1$ on $\partial_{\rm ext}(\Ga)$
    \begin{equation}
    \label{3.1}
W_\Ga(\bar\si_{\partial_{\rm ext}(\Ga)}) \le e^{- c (N_0 \ga^{-1+\alpha +2a} + \ga^{A}|\mathcal S|)}.
     \end{equation}
Same bound holds for minus contours.

\end{thm}

\medskip
We shall prove Theorem \ref{thm3.1} in the rest of the section.  Recall that
the energy of a spin $\si(x,i)$ in the field generated by the configuration $\si'$ outside $(x,i)$
is
\[
H_\ga(\si(x,i)|\si') := -\si(x,i) [h_\ga(x,i;\si') + \ga^{A} (\si(x,i+1)+\si(x,i-1))]
\]
where $h_\ga(x,i;\si'): = \sum_{y\ne x} J_\ga(x,y) \si'(y,i)$.  Then the
Gibbs distribution of   $\si(x,i)$ given $\si'$ is
    \begin{equation}
    \label{3.2}
G_\ga(\si(x,i)|\si')= Z_{\ga,\si'}^{-1} e^{-\beta
H_\ga(\si(x,i)|\si')},
     \end{equation}
with $Z_{\ga,\si'}$ the normalization factor.  The Gibbs conditional
probability of $\si(x,i)$ given $\si'$ and that $\eta(x,i) =1$,
denoted\footnote{We slightly abuse notation here, using the same
symbol for a variable and its possible values.}  by $\mu(\si(x,i)|
\si', \eta(x,i)=1)$, is not always given by \eqref{3.2} because the
condition $\eta(x,i) =1$ involves the spin $\si(x,i)$. However we
obviously have:

\medskip

    \begin{lem}
    \label{lemma3.1}
Let $\si'$ be such that
    \begin{equation}
    \label{3.3}
 |\frac 1{\ell_-}\sum_{y\in C^{\ell_-,i}_x, y\ne x}  \si'(y,i) - m_\beta| < \zeta - \frac 1{\ell_-}.
     \end{equation}
Then
    \begin{equation}
    \label{3.4}
\mu\Big(\si(x,i)| \si', \eta(x,i)=1\Big) = G_\ga(\si(x,i)|\si').
     \end{equation}

    \end{lem}

\medskip
The next lemma gives an upper bound for the probability of violating condition \eqref{3.3}.

\medskip

    \begin{lem}
    \label{lemma3.2}
There are $b<1$ and $c_b>0$ so that for all $\ga$ small enough the
following holds.  Let $\si'$ be a configuration in the complement of
$C^{\ell_-,i}_x$ such that $\eta(y,i)=1$ for all $y$ such that
$(y,i) \notin C^{\ell_-,i}_x$. Denote by $\mu(\,\cdot\, | \;\si',
\eta(x,i)=1)$ the Gibbs conditional probability on
$\{-1,1\}^{C^{\ell_-,i}_x}$ given $\si'$ and that $ \eta(x,i)=1$.
Then (recalling \eqref{2.3} for notation)
    \begin{equation}
    \label{3.5}
\mu\Big(|\si^{(\ell_-)}(x,i)-m_\beta| > b \zeta \;|\; \si', \eta(x,i)=1\Big) \le e^{- c_b \ell_-\zeta^2 }.
     \end{equation}
 \end{lem}

\medskip
\begin{proof}
Since the model is translation invariant, we may take $x=i=0$. Let $\s_y$ stand for $\s(y,0)$, $y\in\Z$,
and let $\s^\pm_y=\s'(y,\pm1)$. We write $\clk=C_{k\ell_-}^{\ell_-,0}$, $k\in\Z$.

The relevant Hamiltonian is then, for $\s=(\s_y)_{y\in\clz}$
\begin{equation}
 \label{ham}
 H(\s)=H_c(\s)+H_b(\s)
\end{equation}
where
\begin{eqnarray}
\label{hams}
 H_c(\s)=-\g^\a\sum_{y\in\clz}\s_y\s_0^{(\ell_-,y)},\quad H_b(\s)=-\sum_{y\in\clz}\s_yh_y,
\end{eqnarray}
\begin{equation}
 \label{ham1}
 \s_k^{(\ell_-,y)}=\frac1{\ell_-}\sum_{z\in\clk}J\left(\g(z-y)\right)\s_z,
\end{equation}
and
\begin{equation}
\label{fields}
h_y=\g^\a\sum_{k\in\Z\atop{k\ne0}}\s_k^{(\ell_-,y)}+\g^{A}(\s_y^++\s_y^-).
\end{equation}

Given the conditions on the boundary and on $J$, it is a straightforward matter to check that
there exists a positive constant $\k$ such that for every $y\in\clz$
\begin{equation}
\label{field-approx}
|h_y-m_\b|\leq\zeta+\k\g^\a.
\end{equation}

The claim of the lemma follows readily from the same bound for the probability of the same event
without the conditioning on $\eta(0,0)=1$ (with a possibly different $c_b$), so we will verify the
latter bound only.

We first dominate in the FKG sense from above and below the model in the volume $\clz$ with the
given boundary conditions by appropriate models without pair couplings within $\clz$, only couplings
to the boundary and extra external magnetic fields, 
so that we will indeed have independent spins in $\clz$ subject to a (uniform) external field appropriately
close to $m_\b$.

For a given constant $M>0$ to be fixed later, let $\mu^\pm$ be the Gibbs measures on spin configurations
in $\clz$ with the following Hamiltonians.
\begin{equation}
 \label{hams-1}
 H^\pm(\s)=-\sum_{y\in\clz}\s_yh^\pm_y,
\end{equation}
where
\begin{equation}
\label{ufields}
h^\pm_y\equiv m_\b\pm[\zeta+(M+\k)\g^\a].
\end{equation}

The result will then follow once we show that
\begin{equation}
 \label{fkg}
 \mu^-(\cdot)\leq\mu(\cdot|\s')\leq\mu^+(\cdot)
\end{equation}
(in the FKG sense), where $\mu$ is the Gibbs measure, and that the bound holds for the probabilities
\begin{eqnarray}
\label{bp}
 \mu^+(\s^{(\ell_-)}>m_\b+b\zeta),\\
\label{bm}
 \mu^-(\s^{(\ell_-)}<m_\b-b\zeta),
\end{eqnarray}
for some $b\in(0,1)$, as soon as $\g$ is close enough to $0$,
where $\s^{(\ell_-)}=\frac1{\ell_-}\sum_{y\in\clz}\s_y$.

An upper bound of the form~(\ref{3.5}) for the expression
in~(\ref{bp}) follows readily from well-known large deviation
bounds, say Bernstein inequality (see e.g. Lemma 1, p. 533 in
\cite{LeCam}), once we notice that under $\mu^+$, the spins in
$\clz$ are iid random variables on $\{-1,+1\}$ with mean
\begin{equation}
 \label{mean}
t_\b(m_\b+\zeta+(M+\k)\g^\a)\leq m_\b + t_\b'(m_\b)[\zeta+(M+\k)\g^\a]\leq m_\b +\tb\zeta,
\end{equation}
where $t_\b:\R\to(-1,1)$ is such that $t_\b(x)=\tanh(\b x)$,
and $\tb<1$ as soon as $\g$ is close enough to $0$, since the derivative of $t_\b$
is less than one on $m_\b$ for $\b>1$. A similar argument establishes a similar
bound for the expression in~(\ref{bm}).

It remains to establish~\eqref{fkg}. We will prove the upper bound.
An argument for the lower bound can be made similarly.

\emph{Proof of the upper bound in~\eqref{fkg}}

We will verify {\em Holley's  condition} (see \cite{Holley}), which in this case reduces to the following bound.
Given $\s,\t\in\{-1,+1\}^{\clz}$
\begin{equation}
 \label{h1}
 -H(\s\wedge\t)-H^+(\s\vee\t)\geq-H(\s)-H^+(\t),
\end{equation}
which in turn reduces to
\begin{eqnarray}\nn
 &\sum_{y\in\clz}(\s\wedge\t)_y\{\g^\a(\s\wedge\t)_0^{(\ell_-,y)}+h_y\}
  +\sum_{y\in\clz}(\s\vee\t)_y\{M\g^\a+\th_y\}&\\
  \label{h2}
 &\geq\sum_{y\in\clz}\s_y\{\g^\a\s_0^{(\ell_-,y)}+h_y\}
  +\sum_{y\in\clz}\t_y\{M\g^\a+\th_y\},&
\end{eqnarray}
where $\th_y\equiv m_\b+\zeta+\k\g^\a$.
We first show that
\begin{equation}
 \label{h3}
 \sum_{y\in\clz}(\s\wedge\t)_yh_y+\sum_{y\in\clz}(\s\vee\t)_y\th_y
 \geq
 \sum_{y\in\clz}\s_yh_y+\sum_{y\in\clz}\t_y\th_y,
\end{equation}
which is equivalent to
\begin{equation}
 \label{h4}
 \sum_{y\in\clz}[(\s\vee\t)_y-\t_y]\th_y
 \geq
 \sum_{y\in\clz}[\s_y-(\s\wedge\t)_y]h_y.
\end{equation}
But $(\s\vee\t)_y-\t_y=\s_y-(\s\wedge\t)_y\geq0$ for all $y$, and~\eqref{h4} follows
from~(\ref{field-approx}), and thence~\eqref{h3} holds. It is enough then to show that
\begin{eqnarray}\nn
 &\sum_{y\in\clz}(\s\wedge\t)_y(\s\wedge\t)_0^{(\ell_-,y)}
  +M\sum_{y\in\clz}(\s\vee\t)_y&\\
  \label{h5}
 &\geq\sum_{y\in\clz}\s_y\s_0^{(\ell_-,y)}
  +M\sum_{y\in\clz}\t_y,&
\end{eqnarray}
which is equivalent to
\begin{eqnarray}\nn
 &M\ell_-\sum_{y\in\clz}[(\s\vee\t)_y-\t_y]=M\ell_-\sum_{y\in\clz}[\s_y-(\s\wedge\t)_y]&\\
  \label{h6}
 &\geq\sum_{y,z\in\clz}J\left(\g(z-y)\right)[\s_y\s_z-(\s\wedge\t)_y(\s\wedge\t)_z].
\end{eqnarray}

Let $\L=\L(\s,\t)=\{x\in\clz:\,\s_x>\t_x\}$, $\tl=|\L|$, and $\L^c=\clz\setminus\L$.
Then the expression on the top of~\eqref{h6} equals $2M\ell_-\tl$ and the one in the
bottom equals
\begin{eqnarray}\nn
&2\sum_{y\in\L,z\in\L^c}J\left(\g(z-y)\right)\s_z+2\sum_{y\in\L^c,z\in\L}J\left(\g(z-y)\right)\s_y&\\
\nn
&\leq
 2\sum_{y\in\L,z\in\L^c}J\left(\g(z-y)\right)+ 2\sum_{y\in\L^c,z\in\L}J\left(\g(z-y)\right)&\\
\label{h6b}
&\leq
4\tilde M \tl(\ell_--\tl)\leq4\tilde M \ell_-\tl,&
\end{eqnarray}
where $\tilde M=\sup_{|r|\leq\g^\a}J(r)$. We conclude that~\eqref{h6} holds as soon as
\begin{equation}
 \label{M}
 M>2J(0)
\end{equation}
and $\g$ is close enough to $0$. Let us then fix an $M$ satisfying~\eqref{M}. We may conclude
that Holley's condition is verified for all $\g$ close enough to $0$, and thence so is
the upper bound in~\eqref{fkg}.
\end{proof}

\medskip

{\bf Remarks.}

\begin{itemize}

\item
Recall that the interaction range is $\ga^{-1}$ so that the condition
$\eta(y,i)=1$ can be required to hold only in the $\ell_-$
intervals on the $i$-th layer which have distance $\le \ga^{-1}$
from $ C^{\ell_-,i}_x$.

\item By the spin flip symmetry Lemma
\ref{lemma3.2} extends to the case where $\eta(y,i)=-1$ with $m_\beta\to -m_\beta$ in \eqref{3.5}.

\item  Suppose that \eqref{3.3} is violated. Then, for $\ga$ small enough,
$|{\si'}^{(\ell_-)}(x,i)-m_\beta| > b \zeta$
no matter what is the value of $\si'(x,i)$.

\end{itemize}

\medskip

We can now start the proof of the Peierls bound which will be achieved after
several manipulations of the partition function in the numerator
of the fraction on the right hand side of \eqref{2.6}.  The first step is to
eliminate some of the vertical interactions in ${\rm sp}(\Ga)$.
Let $S$ be a $+-$ stripe, ${\rm sp}(S)= \{(x,j): x\in I, j=i,i+1\}$.
Denote by $\si'$ a configuration on the complement of ${\rm sp}( S)$.  By the definition of stripes,
$\si'$ is such that
$\eta =1$ on all the $\ell_-$ intervals  on the layer $i+1$ which have distance $\le \ga^{-1}$ from ${\rm sp}( S)$
and $\eta =-1$ on all the $\ell_-$ intervals on the layer $i$
which have distance $\le \ga^{-1}$ from ${\rm sp}( S)$.
We shorthand by $Z_{S,\si'}$ the partition function on   ${\rm sp}( S)$ with boundary conditions
$\si'$ and constraint $\{\eta=\pm 1\}$ on the upper and respectively lower layers of ${\rm sp}( S)$.
We denote by  $Z^0_{S,\si'}$ the same partition function but with the vertical interaction
among the upper and lower  layers of ${\rm sp}( S)$ removed, the vertical interaction
with the complement of ${\rm sp}( S)$ is instead kept.

\medskip

    \begin{prop}
    \label{prop3.1}
There is $c>0$ so that for all $\ga$ small enough
    \begin{equation*}
Z_{S,\si'} \le  e^{- c \ga^{A}|S| }\;Z^0_{S,\si'}.
     \end{equation*}
\end{prop}

\medskip
\begin{proof}
Let $\mu^\eps_{S,\si'}(\cdot)$ be the Gibbs measure where the
vertical interaction in $S$ is $\eps$ instead of $\ga^{A}$, with
$0<\eps\le\ga^{A}$. We have:

\begin{equation}\label{eq:Z-eps}
  \log\frac{Z_{S,\si'}}{Z^0_{S,\si'}}=\sum_{x\in I}\int_0^{\ga^{A}}\mu^\eps_{S,\si'}(\si(x,i)\si(x,i
+1) )d\eps.
\end{equation}

We compute $\mu^\eps_{S,\si'}(\si(x,i)\si(x,i+1) )$ by first conditioning on  $\si''$,
the configuration restricted to $sp(S)\setminus\{(x,i)(x,i+1)\}$:
\begin{eqnarray}
  \mu^\eps_{S,\si'}(\si(x,i)\si(x,i+1)) &=& \mu^\eps_{S,\si'}\left[\mu^\eps_{S,\si'}(\si(x,i)\si(x,i
+1) )\mid\si',\si'', \eta(x,i)=-1,\eta(x,i+1)=+1 \right] \nonumber\\
   &=&  \mu^\eps_{S,\si'}\left[\mathbf 1_{B_x}\mu^\eps_{S,\si'}(\si(x,i)\si(x,i
+1) )\mid\si',\si''\right]+ O(e^{-c\ell_-\zeta^2})
\end{eqnarray}
where:
\begin{equation}\label{def:Ax}
  B_x:=\left\{\si'': |\frac 1{\ell_-}\sum_{y\in C^{\ell_-,i}_x, y\ne x}  \si''(y,i) + m_\beta| <
  \zeta - \frac 1{\ell_-};\,|\frac 1{\ell_-}\sum_{y\in C^{\ell_-,i+1}_x, y\ne x}  \si''(y,i+1) - m_\beta|
  < \zeta - \frac 1{\ell_-}\right\}
\end{equation}
and we have used that $\mu^\eps_{S,\si'}(B^c_x)<
O(e^{-c\ell_-\zeta^2})$ uniformly in $\eps<\ga^{A}$.

It can been seen that on $B_x$,

\begin{eqnarray}
  |\mu^\eps_{S,\si'}(\si(x,i)\si(x,i+1) |\si',\si'')+m^2_\beta| &=&O( \zeta)
\end{eqnarray}
since the vertical interactions in $x$ are uniformly bounded by $\ga^{A}$.

\vskip .5cm
Summing up in $x\in I$ we conclude the statement
\begin{equation}\label{eq:Z-eps-1}
  \log\frac{Z_{S,\si'}}{Z^0_{S,\si'}}\le -|I| \ga^{A} [m_\beta^2 -O(\zeta)].
  \end{equation}
\end{proof}
\vskip .5cm

As an immediate corollary of Proposition \ref{prop3.1} we have:

\medskip
    \begin{cor}
    \label{cor3.1}
Denote by $Z^{0,\mathcal S}_{c(\Ga);\bar\si_{\partial_{\rm ext}(\Ga)}}(\{\si_{c(\Ga)} \Rightarrow \Ga \})$
the partition function  in the numerator
of \eqref{2.6} with the vertical interaction
among the upper and lower  layers of all ${\rm sp}( S)$, $S\in \mathcal S$, removed.  Then
for all $\ga$ small enough
    \begin{equation}
    \label{3.6}
Z_{c(\Ga);\bar\si_{\partial_{\rm ext}(\Ga)}}(\{\si_{c(\Ga)} \Rightarrow \Ga \}) \le
e^{- c \ga^{A}|\mathcal S| }\;Z^{0,\mathcal S}_{c(\Ga);\bar\si_{\partial_{\rm ext}(\Ga)}}(\{\si_{c(\Ga)} \Rightarrow \Ga \})
     \end{equation}
with $c$ as in Proposition \ref{prop3.1}.

    \end{cor}

\medskip

Denote by $Z^{0}_{c(\Ga);\bar\si_{\partial_{\rm ext}(\Ga)}}(\{\si_{c(\Ga)} \Rightarrow \Ga \})$
the partition function  in the numerator
of \eqref{2.6} where it has been removed
the vertical interaction between any two  intervals $C^{\ell_+,i+1}_x$ and
$C^{\ell_+,i}_x$  both in
${\rm sp}(\Ga)$ such that either (i)  $\Theta$ has opposite sign (i.e.\ they belong
to a stripe) or (ii) $\Theta=0$ at least on one of them.

\medskip
    \begin{cor}
    \label{cor3.2}
Let $Z^{0}_{c(\Ga);\bar\si_{\partial_{\rm ext}(\Ga)}}(\{\si_{c(\Ga)} \Rightarrow \Ga \})$
be as above. Then
for all $\ga$ small enough
    \begin{equation}
        \label{3.7}
Z_{c(\Ga);\bar\si_{\partial_{\rm ext}(\Ga)}}(\{\si_{c(\Ga)} \Rightarrow \Ga \}) \le  e^{- c \ga^{A}|\mathcal S|
+ 2\ga^{A}\ell_+ N_0}
\;Z^{0}_{c(\Ga);\bar\si_{\partial_{\rm ext}(\Ga)}}(\{\si_{c(\Ga)} \Rightarrow \Ga \})
     \end{equation}
($c$ as in Proposition \ref{prop3.1} and $N_0$ the number of $\mathcal{D}^{\ell_+}$ intervals in ${\rm sp}(\Ga)$
where $\Theta=0$).

    \end{cor}

   \medskip

\vskip1cm

 Call
    \begin{equation}
    \label{3.8}
\Delta:= \Big\{ (x,i) \in {\rm sp}(\Ga): \Theta(x,i)=0\Big\},\quad |\Delta| = \ell_+N_0
     \end{equation}
and denote by $\si_\Delta$ and $\si'$
the spin configurations in
$\Delta$ and respectively outside $\Delta$.
Since we have dropped all vertical interactions involving spins in $\Delta$
the  system has only
Kac interactions.  A lot is known about such systems and most of what follows
is in fact taken from the existing literature.
We fix $\si'$  outside $\Delta$ and need to bound
   \begin{equation}
    \label{3.9}
Z^0_{\Delta,\si'}( \eta = \eta_\Ga):=  \sum_{\si_{\Delta}} \mathbf 1_{\{\eta = \eta_\Ga\;{\rm on} \Delta\}}
e^{-\beta H_\ga^0(\si_\Delta|\si')}.
     \end{equation}
Observe that     $Z^0_{\Delta,\si'}( \eta = \eta_\Ga)$ factorizes into
a product of partition functions on each layer so that our next estimates will be one-dimensional.

Next step is to coarse-grain to reduce the bound of \eqref{3.9} to
a variational problem involving a free energy functional defined on functions $m(r,i)$,
$r\in \mathbb R, i \in \mathbb Z$.
The scale of the coarse-graining should be chosen to have an error small when compared
to the gain term in \eqref{3.1}: a possible choice that we shall adopt is
$\ell=\ga^{-1/2}$ (which for simplicity we suppose
in $\{2^n, n\in \mathbb N\}$).

As a rule we add a $*$ when we go from the discrete to the continuum, so that
$\Delta^*$ denotes the union
over $(x,i)\in \Delta$ of the
unit intervals $\{(r,i): x\le r< x+1\}$.
We then have (see  Theorem 4.2.2.2 in \cite{presutti})
   \begin{equation}
    \label{3.10}
\log Z^0_{\Delta,\si'}( \eta = \eta_\Ga) \le
-\beta\; \inf_{m_{\Delta}\in \mathcal A}
 F_{\ga,\Delta^*}\left(m_{\Delta}|{\si'}^{(\ga^{-1/2})}\right)+
\beta c \ga^{1/2}\log\ga^{-1}|\Delta|
     \end{equation}
where $m_{\Delta} \in L^\infty(\Delta^*,[-1,1])$; $\si^{(\ga^{-1/2})}$ is
the analogue of $\si^{(\ell_-)}$ in  \eqref{2.3} with $\ell_-$ replaced by
$\ga^{-1/2}$;
$\mathcal A$ is the set of functions $m$ so that for any $(x,i)\in \Delta$ the difference
\[
|\frac 1{\ell_-}\int_{r'}^{r'+\ell_-}  m(r,i) dr \mp m_\beta|,\quad r' = h\ell_- \le x <(h+1)\ell_-
\]
is smaller or larger than $\zeta$ according to the value of  $\eta_\Ga(x,i)$;\\
$\dis{ F_{\ga,\Delta^*}\left(m_{\Delta}|m_{{\Delta}^c}\right) =
 F_{\ga,\Delta^*}(m_{\Delta}) - \sum_i\int
\int  \mathbf 1_{\{(r,i)\in \Delta^*, (r',i)\notin {\Delta^*}\}}
J_\ga(r,r')\,m_{\Delta}(r)m_{{\Delta}^c}(r')}drdr'$, \\
where
       \begin{eqnarray}
    \label{3.11}
&&    F_{\ga,\Delta^*}(m_{\Delta})=  -\frac 12 \sum_i\int_{\{(r,i)\in \Delta^*\}}
\int_{\{(r',i)\in \Delta^*\}}J_\ga(r,r')\,m_{\Delta} (r)m_{\Delta} (r')dr dr' \nonumber 
- \frac 1\beta  \int_{\Delta^*} I(m_{\Delta}(r))dr
\\&&
I(m) = - \frac{1-m}{2} \log\,\frac{1-m}{2} - \frac{1+m}{2}
\log\,\frac{1+m}{2};
      \end{eqnarray}
 finally $c$ in \eqref{3.10}
is a constant.

Observe that the last term in \eqref{3.10} is bounded by
$\beta c N_0 \ga^{-1/2-\alpha}\log\ga^{-1}$,
thus the ``error'' in \eqref{3.10} is ``small'' with respect to the gain term in \eqref{3.1}
(because $a$ and $\alpha$ are suitably small).

The next step exploits the stability property of the functional
in a neighborhood of the stationary profiles identically equal to $m_\beta$
(or to $-m_\beta$). The intersection of a layer   $\{(r,i): r \in \mathbb R\}$
with $\Delta^*$ (supposing it is non empty) is made of consecutive disconnected
intervals
\[
I_{h,i} = [(r,i): r'_h \le r <r''_h),\quad r'_h, r''_h \in \ell_+ \mathbb Z,
\]
where the extremes of the separating intervals $[r''_h,r'_{h+1})$  are either the endpoints of
a stripe layer, or the intersection with $\{(r,i): r \in \mathbb R\}$ of an interior
${\rm int}_j(\Ga) ^*$.  Thus
by construction $\theta (r'_h,i) = \pm 1$ and
since $m_\Delta \in \mathcal A$,
for all $k$ such that
$ [k\ell_-,(k+1)\ell_-)\subseteq [r'_h, r'_h+\ell_+)$ we have either
\[
\Big|\frac 1{\ell_-} \int_{k\ell_-}^{(k+1)\ell_-} m_{\Delta}(r) dr - m_\beta\Big| \le \zeta
\]
or the same with $m_\beta\to -m_\beta$.  The analogous property holds in
$[r''_h-\ell_+,r''_h)$.  Let us focus for instance on the interval
$[r'_h,r'_h+\ell_+)$, call $r_{{\rm mid}}:=  r'_h+\ell_+/2$ and, to fix the ideas,  suppose the averages of
$m_\Delta$ are close to $m_\beta$. Then by Theorem 6.3.3.1 in \cite{presutti} there are $\om>0$ and $c$ so that
the inf in \eqref{3.10} is achieved
on functions $m$ with the following property.
  \begin{equation}
    \label{3.12}
\sup_{|r-r_{{\rm mid}}| \le \ga^{-1}} |m(r)-m_\beta| \le c e^{-{\om}  \ga \ell_+} = c e^{- {\om} \ga^{-\alpha}}.
     \end{equation}
Thus,
   \begin{eqnarray*}
\dis{\inf_{m_{\Delta}\in \mathcal A}
 F_{\ga,\Delta^*}\left(m_{\Delta}|{\si'}^{(\ga^{-1/2})}\right) }
\ge&&
 \inf_{m_{\Delta}\in \mathcal A; m_\Delta= m_\beta\; {\rm on}\;
 |r-r_{{\rm mid}}| \le \ga^{-1}}
 F_{\ga,\Delta^*}\left(m_{\Delta}|{\si'}^{(\ga^{-1/2})}\right)\\
 &-&\ga^{-1}c' e^{- {\om} \ga^{-\alpha}}.
     \end{eqnarray*}
By changing the constant
$c$ in \eqref{3.10} we can then restrict  in \eqref{3.10}
to functions which are identically equal to
$m_\beta$ or to $-m_\beta$ depending on the
value of $\eta_\Ga$ in all the intervals of the form
$|r-r_{\rm mid}| \le \ga^{-1}$ with
$r_{\rm mid}$ at distance $\ell_+/2$ from an endpoint of
any of the
$I_{h,i}$.

Call $I_1= [r'_h,r'_h+\frac{\ell_+}2]$,
$I_2=[r''_h-\frac{\ell_+}2, r''_h]$ and $I_0=I_{h,i}\setminus\{I_1\cup
I_2\}$. Let $m$ be a function on $I_{h,i}$  equal to $\pm m_\beta$
in the two intervals $|r-r_{\rm mid}| \le \ga^{-1}$, with
$r_{\rm mid}$ at distance $\ell_+/2$ from $r'_h$ and from $r''_h$, respectively.
Call $m_1$, $m_2$ and $m_0$ the restriction of $m$ to $I_1$, $I_2$ and $I_0$.  We then have
       \begin{eqnarray*}
&&    F_{\ga,I_{h,i}}(m|{\si'}^{(\ga^{-1/2})})=  F_{\ga,I_{1}}(m_1|{\si'}^{(\ga^{-1/2})})
+F_{\ga,I_{2}}(m_2|{\si'}^{(\ga^{-1/2})}) +\mathcal F_{\ga,I_0}(m_0) - 2C m_\beta^2
\end{eqnarray*}
where
\begin{eqnarray}
&&
C = \frac 12 \int_{I_0} \int _{I_1} J_\ga(r,r')drdr'= \frac 12 \int_{I_0} \int _{I_2} J_\ga(r,r')drdr',\quad
f_\beta(m) = - \frac{m^2}2 - \frac{1}{\beta} I(m),\nonumber\\
&&\mathcal F_{\ga,I_0}(m_0) = \int_{I_0} f_\beta (m_0(r))dr + \frac{\beta}4 \int_{I_0}\int_{I_0}
J_{\ga}(r,r') \Big( m_0(r)-m_0(r')\Big)^2drdr'.
      \end{eqnarray}

By Theorem 6.4.2.3 in \cite{presutti}
  \begin{equation}
    \label{3.14}
\mathcal F_{\ga,I_0}(m_0) \ge |I_0| f_\beta(m_\beta) + c \ell_- \zeta^2 (2n+p),
     \end{equation}
where $p$ is the number of intervals $C^{\ell_-,i}\subset I_0$ where $\eta_\Ga =0$
and $n$ is the number of consecutive pairs of intervals in $I_0$ where $\eta_\Ga$ changes
from $1$ to $-1$ or viceversa.     We can then rewrite
  \begin{equation}
    \label{3.15}
\mathcal F_{\ga,I_0}(m_0) \ge \mathcal F_{\ga,I_0}(m_\beta \mathbf 1_{I_0}) + c \ell_- \zeta^2 (2n+p).
     \end{equation}
Call $\tilde m_1= m_1$ if $\eta_\Ga=1$ on $I_1$ and $=-m_1$ otherwise,
analogous notation are used for $m_2$; similarly call
$\si''$ the configuration outside $\Delta$ obtained from
$ \si'$ by flipping the spins in ${\rm int}_k^-$ and in the parts of
the stripes where $\Theta=-1$.
Then
calling $\tilde m$
the function equal to $m_\beta$ on $I_0$ and to $\tilde m_1$ and $\tilde m_2$ on $I_1$ and $  I_2$
  \begin{equation}
    \label{3.16}
 F_{\ga,I_{h,i}}(m|{\si'}^{(\ga^{-1/2})}) \ge F_{\ga,I_{h,i}}(\tilde m|{\si''}^{(\ga^{-1/2})})
 +
  c \ell_- \zeta^2 (2n+p).
     \end{equation}
By collecting the above bounds on all the intervals $I_{h,i}$ we then get from \eqref{3.10}
   \begin{equation}
    \label{3.17}
\log Z^0_{\Delta,\si'}( \eta = \eta_\Ga) \le
-\beta\;
 F_{\ga,\Delta^*}\left(\tilde m_{\Delta}|{\si''}^{(\ga^{-1/2})}\right)+
\beta c \ga^{1/2}\log\ga^{-1}|\Delta| -  c^{*} \ell_- \zeta^2 \frac{|\Delta|}{\ell_+},
     \end{equation}
where $\tilde m_{\Delta}$ is such that its $\ell_-$ averages are all close to $m_\beta$,
$\si''$ is obtained from $\si'$ by flipping the spins in all minus interiors of ${\rm sp}(\Ga)$
and in the minus parts of the stripes; instead
$\si''=\si'$ in the plus interiors and in the plus parts of the stripes. Finally
the sum of the numbers $(2n+p)$ over all the intervals $I_{h,i}$ is bounded proportionally by a
factor $1/K$ to the number of $C^{\ell_+,i}$ intervals in $\Delta$,
and $c^{*}= c/K$.

Using again Theorem 4.2.2.2 in \cite{presutti}, we have
     \begin{equation}
    \label{3.18}
    \log Z^0_{\Delta,\si''}( \eta = 1) \ge
-\beta\;  F_{\ga,\Delta^*}\left(\tilde m_{\Delta}|{\si''}^{(\ga^{-1/2})}\right) -\beta c
 \ga^{1/2}\log \ga^{-1}|\Delta|
     \end{equation}
so that for $\ga$ small enough
   \begin{equation}
    \label{3.19}
 Z^0_{\Delta,\si'}( \eta = \eta_\Ga) \le      Z^0_{\Delta,\si''}( \eta = 1)\times
e^{- \frac{ c^{*}}2 \ga^{-1+ \alpha + 2a} N_0 }.
     \end{equation}
We have thus proved that $Z_{c(\Ga);\bar\si_{\partial_{\rm ext}(\Ga)}}(\{\si_{c(\Ga)}
\Rightarrow \Ga \})$ is bounded by
         \begin{equation}
    \label{3.20}
 Z^0_{c(\Ga);\bar\si_{\partial_{\rm ext}(\Ga)}}(\{\si_{c(\Ga)}: \Theta = 1\; \rm{on}\;
{\rm sp}(\Ga)\; {and \; all}\, \partial_{k}^{\pm}(\Ga)\})\times e^{- c
\ga^{A}|\mathcal S| + 2\ga^{A}\ell_+ N_0} e^{- \frac{ c^{*}}2
\ga^{-1+ \alpha + 2a} N_0 },
     \end{equation}
where the superscript $ Z^0$ recalls that in the partition function
some vertical interactions are missing: the missing ones are those
between the layers of the stripes $S\in \mathcal S$ and those
involving the $(x,i)\in {\rm sp}(\Ga)$ where $\Theta=0$.  A proof
analogous to that of Proposition \ref{prop3.1} shows that if $I=
[k\ell_+,(k+1)\ell_+)$, $S = \{(x,j): x\in I, j\in \{i,i+1\}$,
$\si'$ a spin configuration outside $S$ with $\Theta \equiv 1$ and
$Z^0_{S,\si'} $ the partition function in $S$ with the constraint
$\theta=1$ identically and without vertical interaction, then there
is $c>0$ so that for all $\ga$ small enough
    \begin{equation*}
Z^0_{S,\si'} \le  e^{- c \ga^{A}|S| }\;Z_{S,\si'}
     \end{equation*}
where in the latter the vertical interaction is present.  Applying repeatedly this inequality
we then get from \eqref{3.20} the proof of Theorem \ref{thm3.1}.
\vskip1cm

\setcounter{equation}{0}

\section{Peierls estimates}

In this section we prove the following theorem from which \eqref{1.5} follows at once for $\ga$ small enough.

\begin{thm}
\label{thm:PE}
In the notation of Theorem \ref{thm3.1},  a positive constant $\tilde c$ can be found so that for all $\ga$ small,
\begin{equation}\label{PE}
  \sum_{\Ga: \rm{sp}(\Ga)\ni 0}W_\Ga(\bar\si_{\partial_{\rm ext}(\Ga)})< e^{- \tilde c \ga^{-1+\alpha +2a} },
\end{equation}
where $\alpha$ and $a$ are the same as in Theorem \ref{thm3.1}.
\end{thm}

\vskip .5cm


\begin{proof}
In the notation of Theorem \ref{thm3.1}, if $\Gamma$ is a plus contour
we may rewrite \eqref{3.1} as follows
$$
 \label{mumu}
 W_{\Gamma}(\bar\si_{\partial_{\rm ext}(\Ga)}) \le \prod_{I \in \mathcal{I}_0} e^{-\frac{c}{2}
 \gamma^{-1+\alpha+2a}} \prod_{S\in \mathcal S}e^{-\frac{c}{4} \gamma^{-1+\alpha+2a}-c\gamma^{A}|S|},
$$
where $\mathcal{I}_0$ is the set of $\mathcal{D}^{\ell_+}$ intervals in $\rm{sp}(\Gamma)$ with $\Theta=0$, and
we have used that next to each side of $S\in \mathcal S$, and in at least one of the layers, there must be an
interval in $\mathcal I_0$. Thus a simple correspondence can be established in such a way that each such
interval is ``used" by at most 2 stripes in $\mathcal S$.

Our goal is to show that for suitable $\psi >0$ small (see \eqref{mumu4}) and all $\gamma$
small
\begin{equation}
\sum_{\Gamma \colon 0 \in \rm{sp}(\Gamma)} W_{\Gamma}(\bar\si_{\partial_{\rm ext}(\Ga)}) <\psi.
\label{mumu1}
\end{equation}
The sum over all $\Gamma$ so that $0 \in \rm{sp}(\Gamma)$ can be obtained by summing over trees
where each vertex in the tree corresponds to an $I \in \mathcal I_0$ or to $\rm{sp}(S)$ for $S \in \mathcal S$, and which
will cover $\rm{sp}(\Gamma)$ exactly (a spanning tree); the types depend also on $\eta_{\Gamma}$.
At each step, the number of descendants in the next generation is bounded by the number of connected
sites in $\mathbb{Z}\times\mathbb{Z}$, i.e. at most $8$ in case of an $I$, and at most $2|S|+8$
in case of an $S$. We may span the tree from a root, and each next generation of a vertex is formed by
vertices in correspondence to connected $I$ or $S$ in $\rm{sp}(\Gamma)$ that have not yet appeared.

The root can be thought to be the $I$ or $S$ that contains the
origin. For an $I$ we use the crude bound $3^{\ell_+/\ell_-}$ for
the number of possibilities with $\Theta=0$ (taking all
possibilities for the $\eta$ variables). For an $S$ the number of
possibilities is at most $4|S|$ (by considering the location of the
origin in ${\rm{sp}}(S)$ and the type of $S$). To achieve
\eqref{mumu1}, it suffices to have for such a small positive $\psi$;
\begin{equation}
(1+\psi)^{8}  e^{-\frac{c}{2}\gamma^{-1+\alpha
+2a}}3^{\ell_+/\ell_-}+ \sum_{S\colon 0 \in {\rm{sp}}(S)}
(1+\psi)^{2|S|+8} e^{-c\gamma^{A}|S|} e^{-\frac{c}{4}
\gamma^{-1+\alpha+2a} } <\psi. \label{mumu2}
\end{equation}
Indeed, for \eqref{mumu1} it suffices to prove that the sum for all
trees with at most $m$ generations is bounded by $\psi$, for all
$m$. This is done by induction on $m$. We can see it at once
 by treating the simple cases the trees are only the root ($m=0$) or have one generation, and
then by expanding depending on the first generation. Indeed, when $m=0$ the tree is only the root
and the bound becomes
\begin{equation*}
e^{-\frac{c}{2}\gamma^{-1+\alpha +2a}}3^{\ell_+/\ell_-}+
\sum_{S\colon 0 \in {\rm{sp}}(S)} e^{-c\gamma^{A}|S|}
e^{-\frac{c}{4} \gamma^{-1+\alpha+2a} },
\end{equation*}
which would be bounded by $\psi$.
Upon conditioning on the first generation and using that the sum starting on each such nodes
is bounded by $\psi$ (by the induction assumption), the induction follows easily.
This is the reason for the factors $(1+\psi)^{8}$ in case of an $I$ or
$(1+ \psi)^{2|S|+8}$ in case of an $S$.

It remains  to check the validity of \eqref{mumu2}. We can see it by breaking into two:

\begin{equation}
(1+\psi)^{8} e^{-\frac{c}{2}\gamma^{-1+\alpha +2a}}3^{\ell_+/\ell_-} < \psi /2
\label{mumu3}
\end{equation}
and
\begin{equation}
\sum_{{S: 0 \in {\rm{sp}}(S)}} (1+\psi)^{2|S|+8} e^{-c\gamma^{A}|S|}
e^{-\frac{c}{4} \gamma^{-1+\alpha+2a} } <{\psi}/{2}.
\end{equation}

Since we assumed that  $\alpha$ and $a$ are suitably small, we easily see that the first estimate is achieved
(for all $\gamma$ small) by taking $\psi$ of the order $e^{-\tilde c \gamma^{-1+\alpha+2a}}$ for
$\tilde c < c/4$. For the second one needs to see
$$
\sum_{n \ge 1} 4n (1+\psi)^{2n+8} e^{-nc\gamma^{A}} e^{-\frac{
c}{4}\gamma^{-1+\alpha +2a}}<\psi/2,
$$
which boils down to show that
$$
8(1+\psi)^{10}e^{-\frac{c}{4}\gamma^{-1+\alpha+2a}}e^{-c\gamma^{A}}<
\psi (1-(1+\psi)^2 e^{-\gamma^{A}})^2
$$
and we can check that both work for
\begin{equation}
\label{mumu4}
\psi =e^{-\tilde c \gamma^{-1+\alpha+2a}}
\end{equation}
with suitable $\tilde c>0$.
\end{proof}

\setcounter{equation}{0}

\section{Proof of Theorem \ref{thm1.1}}
\label{sec:5}

Let $\La_n$ be any increasing sequence of $\mathcal
D^{\ell_+}$-measurable regions invading $\mathbb Z\times \mathbb Z$
and let $\mu_{\ga,\La_n;\bar\si_{\La_n^c}}^{\pm}$ be Gibbs measures
with boundary conditions $\bar\si_{\La_n^c}$ such that $\Theta$ is
identically 1 (respectively $-1$) on the complement $\La_n^c$ of
$\La_n$. By general arguments based on the validity of the Peierls
bounds, see \cite{BMP} and Chapter 12 in \cite{presutti},
$\mu_{\ga,\La_n;\bar\si_{\La_n^c}}^{\pm}$ converge weakly,
independently of the choice of $\La_n$ and of the boundary
conditions, to distinct DLR measures that we denote by
$\mu_{\ga}^{\pm}$ (the statement would follow from ferromagnetic
inequalities if the plus/minus boundary conditions were realized by
spin configurations identically equal to 1, respectively $-1$).  By
the arbitrariness of the sequence $\La_n$ and of the boundary
conditions it then follows that $\mu_{\ga}^{\pm}$ are invariant
under horizontal translations by multiples of $\ell_+$ and under
vertical translations.  As a consequence any translational invariant
DLR measure $\mu$ is a convex combination of $\mu_{\ga}^{\pm}$: this
is based on an extension of the original proof by Gallavotti and
Miracle-Sole for the Ising model at small temperatures, see again
\cite{BMP} and Chapter 12 in \cite{presutti}.

Since any weak limit $\mu$ of $\mu_{\ga,\La}^{\rm per}$ is invariant under translation,
then $\mu = a\mu_{\ga}^{+}+(1-a)\mu_{\ga}^{-}$; by the spin flip symmetry
$\mu(\si(0,0)=1)= \frac 12$ hence $a= \frac 12$ and Theorem \ref{thm1.1} is proved.

\vskip 1.5cm

\noindent {\bf Acknowledgement}

MEV thanks the warm hospitality of GSSI, L'Aquila, where
part of this research was done.

Research partially supported by CNPq grant 474233/2012-0.
MEV's work is partially supported by CNPq grant 304217/2011-5 and
Faperj grant E-24/2013-132035.
LRF's work is partially supported by CNPq grant 305760/2010-6 and
Fapesp grant 2009/52379-8.

\end{document}